\documentclass[11pt]{article}

\usepackage{tgpagella}
\usepackage[utf8]{inputenc}
\usepackage[T1]{fontenc}
\usepackage{amssymb}
\usepackage{amsfonts}
\usepackage{amsmath}
\usepackage{amsthm}
\usepackage{mathrsfs}
\usepackage{comment}

\newtheorem{theorem}{Theorem}
\newtheorem{proposition}[theorem]{Proposition}

\newtheorem{cor}[theorem]{Corollary}

\theoremstyle{definition}
\newtheorem{definition}[theorem]{Definition}
\newtheorem{convention}[theorem]{Convention}

\newcommand{\PA}{\textnormal{PA}}

\newcommand{\ITB}{\textnormal{ITB}}
\newcommand{\set}[2]{\left\{ #1 \ \mid \ #2 \right\} }

\newcommand{\CT}{\textnormal{CT}}
\newcommand{\df}[1]{\textbf{#1}}
\newcommand{\num}[1]{\underline{#1}}

\newcommand{\INT}{\textnormal{INT}}

\newcommand{\ElDiag}{\textnormal{ElDiag}}
\newcommand{\val}[1]{{#1}^{\circ}}

\renewcommand{\Pr}{\textnormal{Pr}}

\newcommand{\QFC}{\textnormal{QFC}}
\newcommand{\DC}{\textnormal{DC}}
\newcommand{\ACDC}{\textnormal{ACDC}}

\newcommand{\qfSent}{\textnormal{qfSent}}
\newcommand{\Tr}{\textnormal{Tr}}

\newcommand{\Prop}{\textnormal{Prop}}

\newcommand{\qcr}[1]{\ulcorner #1 \urcorner}

\newcommand{\LPA}{\mathscr{L}_{\PA}}
\newcommand{\form}{\textnormal{Form}}
\newcommand{\Term}{\textnormal{Term}}
\newcommand{\Sent}{\textnormal{Sent}}
\newcommand{\Var}{\textnormal{Var}}
\newcommand{\ClTerm}{\textnormal{ClTerm}}
\newcommand{\TermSeq}{\textnormal{TermSeq}}
\newcommand{\ClTermSeq}{\textnormal{ClTermSeq}}
\newcommand{\SentSeq}{\textnormal{SentSeq}}
\newcommand{\Num}{\textnormal{Num}}
\newcommand{\FV}{\textnormal{FV}}

\newcommand{\GRP}{\textnormal{GRP}}
\newcommand{\PC}{\textnormal{PC}}
\newcommand{\PS}{\textnormal{PS}}
\newcommand{\IDelta}{\textnormal{I}\Delta}

\newcommand{\Unique}{\bigcirc}
\newcommand{\Comp}{\textnormal{Comp}}
\newcommand{\Asn}{\textnormal{Asn}}

\title{Compositional truth with propositional tautologies and quantifier-free correctness}
\author{Bartosz Wcisło}

\begin{document}

	\maketitle
	
\begin{abstract}
	In \cite{cieslinskict0}, Cieśliński asked whether compositional truth theory with the additional axiom that all propositional tautologies are true is conservative over Peano Arithmetic. We provide a partial answer to this question, showing that if we additionally assume that truth predicate agrees with arithmetical truth on quantifier-free sentences, the resulting theory is as strong as $\Delta_0$-induction for the compositional truth predicate, hence non-conservative. On the other hand, it can be shown with a routine argument that the principle of quantifier-free correctness is itself conservative. 
\end{abstract}

\section{Introduction}

It is a very widespread phenomenon in logic that if a theory $S_1$ can formulate a truth predicate for a theory $S_2$, then $S_1$ is stronger than $S_2$, a claim which can be made precise in many different ways.

This phenomenon, stripped down to its essence, is investigated in the area of truth theory. Truth theories are axiomatic theories which arise by adding a fresh predicate $T(x)$ to a base theory $B$ which handles syntactic notions (Peano arithmetic, $\PA$, is an example of such a theory). The intended interpretation of $T$ is the set of (codes of) true sentences of the base theory. By considering various possible axioms governing the behaviour of $T$, we investigate the impact of various notions of truth on the properties on the obtained theory.

One line of research in this area asks what precise properties of the truth predicate make a theory with a truth predicate non-conservative over the base theory. (A theory $S_1$ is conservative over its subtheory $S_2$ if it does not prove any theorems in the language of $S_2$ which are not already provable in that subtheory.) 

It is rather straightforward to see that if we add to $\PA$ a truth predicate which satisfies compositional axioms and the full induction scheme in the arithmetical language extended with the unary truth predicate, then by induction on lengths of proofs we can show that all theorems of $\PA$ are true and hence arithmetic is consistent. On the other hand, by a nontrivial result of Kotlarski, Krajewski, and Lachlan from \cite{kkl}, the theory of pure compositional truth predicate with no induction is conservative over $\PA$. Recent research brought much better understanding of which exact principles weaker than full induction yield a nonconservative extension of $\PA$.\footnote{A comprehensive discussion of recent discoveries can be found in \cite{cies_ksiazka}.}

One of the persistent open questions in this line of research asks whether compositional truth theory over $\PA$ with an additional axiom expressing that all propositional tautologies are true is conservative over arithmetic. We know that related principles such as "truth is closed under propositional logic" or "valid sentences of first-order logic are true" are not conservative and indeed are all equivalent to $\Delta_0$-induction for the truth predicate.\footnote{The question was originally stated by Cieśliński in \cite{cieslinskict0}. It was also asked by Enayat and Pakhomov in \cite{EnayatPakhomov}.} 

In this article, we provide a partial answer to Cieśliński's question. We show that $\CT^-$ extended with the principle expressing that propositional tautologies are true becomes nonconservative upon adding quantifier-free correctness principle $\QFC$ which states that $T$ predicate agrees with partial arithmetical truth predicates on quantifier-free sentences. The principle $\QFC$ can itself be easily seen to be conservative over $\PA$ (we include a proof in the Appendix B; it is routine).

Our result can therefore be seen as a certain no-go theorem. Our methods for showing conservativity of truth theories behave very well when we demand that several such properties are satisfied at once. Therefore our theorem seems to impose certain restriction on what methods can be used to attack the problem of propositional tautologies.

\section{Preliminaries}

\subsection{Arithmetic}
In this paper, we consider truth theories over Peano Arithmetic ($\PA$) formulated in the language $\{+,\times, S,0\}$. It is well known that $\PA$, as well as its much weaker subsystems, are capable of formalising syntax. This topic is standard and the reader can find its discussion e.g. in \cite{kaye} or \cite{hajekpudlak}. Below, we list some formulae defining formalised syntactic notions which we will use throughout the paper. 

\begin{definition} \label{def_formalised_syntax} \
	\begin{itemize}
		\item $\Var(x)$ defines the set of  (codes of) first-order variables.
		\item $\Term_{\LPA}(x)$ defines the set of (codes of) terms of the arithmetical language.
		\item $\ClTerm_{\LPA}(x)$ defines the set of (codes of) closed terms of the arithmetical language. 
		\item $\Num(x,y)$ means that $y$ is (the code of) the canonical numeral denoting $x$. We will use the expression $y = \num{x}$ interchangeably.
		\item $\val{t}=x$ means that $t$ is (a code of) a closed arithmetical term and its formally computed value is $x$. 
		\item $\form_{\LPA}(x)$ defines the set of (codes of) arithmetical formulae.
		\item $\form_{\LPA}^{\leq 1}(x)$ defines the set of (codes of) arithmetical formulae with at most one free variable.
		\item $\Sent_{\LPA}(x)$ defines the set of (codes of) arithmetical sentences. 
		\item $\SentSeq_{\LPA}(x)$ defines the set of (codes of) sequences of arithmetical sentences.
		\item $\qfSent_{\LPA}(x)$ defines the set of (codes of) quantifier-free arithmetical sentences. 
		\item $\Pr_{\PA}(d,\phi)$ means that $d$ is (a G\"odel code of) a proof of $\phi$ in $\PA$. $\Pr_{\PA}(\phi)$ means that $\phi$ is provable in $\PA$. 
		\item $\FV(x,y)$ means that $y$ is (a code of) an arithmetical formula and $x$ is amongst its free variables.  
		\item $\Asn(\alpha,x)$ means that $\phi$ is (a code of) an arithmetical term or formula and $\alpha$ is an \df{assignment} for $x$, i.e., a function whose domain contains its free variables.
		\item If $t \in \Term_{\LPA}$ and $\alpha$ is an assignment for $t$, then by $t^{\alpha}=x$, we mean that $x$ is the formally computed value of the term $t$ under the assignment $\alpha$.  
		\end{itemize}
\end{definition}

In the paper, we will make an extensive use of a number of conventions. 

\begin{convention} \label{conv_definicje_syntaktyczne} \
	\begin{itemize}
		\item We will use formulae defining syntactic objects as if they were denoting the defined sets. For instance, we will write $x \in \Sent_{\LPA}$ interchangeably with $\Sent_{\LPA}(x)$.
		\item We will often omit expressions defining syntactic operations and simply write the results of these operations in their stead. For example, we will write $T(\phi \wedge \psi)$ meaning "$\eta$ is the conjunction of (the codes of) the  sentences $\phi, \psi,$ and $T(\eta)$."
		\item We will use formulae defining functions as if they actually were function symbols, e.g. writing $\num{x}$ or $\val{t}$ like stand-alone expressions.
		\item We will in general omit Quine corners and conflate formulae with their G\"odel codes. This should not lead to any confusion.
		\item We will use expressions $x \in \FV(\phi)$ and $\alpha \in \Asn(\phi)$ interchangeably with $\FV(x,\phi)$ and $\Asn(\alpha,\phi)$. Moreover, we will use the expressions $\FV(\phi), \Asn(\phi)$ as if they had a stand-alone meaning, denoting sets of free variables and of $\phi$-assignments respectively.
	\end{itemize}
\end{convention}

In this paper, we analyse compositional truth theory. Let us define the theory in question.

\begin{definition} \label{def_ctminus}
	By $\CT^-$ we mean a theory formulated in the arithmetical language extended with a fresh unary predicate $T(x)$ obtained by adding to $\PA$ the following axioms:
	\begin{enumerate}
		\item $\forall s,t \in \ClTerm_{\LPA} \ \Big(T(s=t) \equiv \val{s} = \val{t}\Big).$
		\item $\forall \phi \in \Sent_{\LPA} \ \Big(T \neg \phi \equiv \neg T \phi \Big).$
		\item $\forall \phi, \psi \in \Sent_{\LPA} \ \Big(T (\phi \vee \psi) \equiv T\phi \vee T \psi \Big).$
		\item $\forall \phi \in \form^{\leq 1}_{\LPA} \forall v \in \FV(\phi) \ \Big(T\exists v \phi \equiv \exists x T\phi(\num{x}) \Big).$
		\item $\forall \bar{s}, \bar{t} \in \ClTermSeq_{\LPA} \forall \phi \in \form_{\LPA} \ \Big( \phi(\bar{s}), \phi(\bar{t}) \in \Sent_{\LPA} \wedge \bar{\val{s}} = \bar{\val{t}} \rightarrow T\phi(\bar{s}) \equiv T \phi(\bar{t})\Big).$
	\end{enumerate}
\end{definition}
Notice that in the axioms of $\CT^-$ we do not assume any induction for the formulae containing the compositional truth predicate. 

\begin{definition} \label{def_ct_ctn}
	By $\CT$ we mean the theory obtained by adding to $\CT^-$ full induction scheme for formulae in the full language (i.e., arithmetical language extended with the unary truth predicate).
	
	By $\CT_n$ we mean $\CT^-$ with $\Sigma_n$-induction in the extended language, for $n \geq 0$.
\end{definition}

It is very well known that $\PA$ (and,in fact, its much weaker fragments) can define partial truth predicates, i.e., formulae which satisfy axioms of $\CT^-$ for sentences of some specific syntactic shape.\footnote{See Chapter I, Section 2(c) of \cite{hajekpudlak}.} In this paper, we will only need a very special case of this fact.

\begin{proposition} \label{stw_qf_partial_truth}
	There exists an arithmetical formula $\Tr_0(x)$ which satisfies axioms 1--3 of $\CT^-$ restricted to $\phi, \psi \in \qfSent_{\LPA}$, provably in $\PA$. 
\end{proposition}

\subsection{The Tarski boundary}

Recall that a theory $S_1$ is \df{conservative} over $S_2$ if $S_1 \supseteq S_2$ and whenever $\phi$ is a sentence from the language of $S_2$ and $S_1 \vdash \phi$, then $S_2 \vdash \phi$. It is a persistent phenomenon in logic that the presence of a truth predicate adds substantial strength to theories in question, as witnessed by the following classical theorem:
\begin{theorem} \label{tw_ct_not_conservative}
	$\CT$ is not conservative over $\PA$.
\end{theorem}
The compositional truth predicate can be employed to prove by induction on the size of proofs that whatever is provable in $\PA$ is true. This allows us to derive the consistency statement for $\PA$ which is unprovable in Peano Arithmetic itself by G\"odel's Second Theorem. The straightforward argument mentioned above uses $\Pi_1$-induction for the compositional truth predicate, but in fact one can do better: 
\begin{theorem} \label{tw_ct0_not_conservative}
	$\CT_0$ is not conservative over $\PA$. 
\end{theorem}
As a matter of fact, as shown in \cite{lelyk_thesis}, $\Delta_0$-induction is equivalent over $\CT^-$ to the following \df{Global Reflection Principle} ($\GRP$):
\begin{displaymath}
\forall \phi \in \Sent_{\LPA} \ \Big(\Pr(\phi) \rightarrow T \phi \Big).
\end{displaymath}
Note that $\GRP$ is, in a way, the exact reason why $\CT$ is not conservative over $\PA$. On the other hand, one of the most important features of $\CT^-$ is that it cannot prove any new arithmetical theorems.
\begin{theorem}[Essentially Kotlarski--Krajewski--Lachlan] \label{tw_kkl}
	$\CT^-$ is conservative over $\PA$.
\end{theorem}
Now, as we can see, compositional truth by itself can be deemed "weak," but it becomes strong upon adding some induction. One of the main goals of our research is to understand what principles can be added to $\CT^-$ in order to make it nonconservative. It turns out that $\CT_0$ plays a crucial role in this research. A number of apparently very distinct principles turns out to be exactly equivalent with $\Delta_0$-induction for the truth predicate. Let us present the one which largely motivates the research in this paper.

\begin{definition} \label{def_propositional_closure}
	By \df{Propositional Closure Principle} $(\PC)$ we mean the following axiom:
	\begin{displaymath}
	\forall \phi \in \Sent_{\LPA} \Big(\Pr^{\Prop}_{T}(\phi) \rightarrow T\phi \Big).
	\end{displaymath}
\end{definition}
The formula $\Pr^{\Prop}_T(x)$ means that $x$ is provable from true premises in propositional logic. 

It was proved in \cite{cies} 
that $\PC$ is actually equivalent over $\CT^-$ to $\CT_0$. This is a very surprising result: the mere closure of truth under propositional logic is actually enough to show that consequences of $\PA$ are true. 

We can form principles similar to $\PC$ which employ stronger closure conditions:
\begin{itemize}
	\item "Truth is closed under provability in propositional logic".
	\item "Truth is closed under provability in first-order logic".
	\item "Truth is closed under provability in $\PA$". 
\end{itemize}
We can also weaken these principles so that they only express soundness of discussed systems, not closure properties.
\begin{itemize}
	\item "Any sentences provable in first-order logic is true".
	\item "Any sentence provable in $\PA$ is true". 
\end{itemize}
It turns out that all the principles listed above are equivalent to each other over $\CT^-$.\footnote{See \cite{cies_ksiazka}.} One axiom which is noticeably absent from the list is the soundness counterpart of $\PC$. This is not an accident. Whether this principle is conservative over $\PA$ is still an open problem. Let us state our official definition.

\begin{definition} \label{def_ps}
	By \df{propositional soundness principle} ($\PS$), we mean the following axiom:
	\begin{displaymath}
	\forall \phi \in \Sent_{\LPA} \Big(\Pr^{\Prop}_{\emptyset}(\phi) \rightarrow T\phi\Big).
	\end{displaymath}
\end{definition}
The formula $\Pr^{\Prop}_{\emptyset}(\phi)$ expresses that $\phi$ is provable in propositional logic from the empty set of premises. In other words, $\PS$ states that any propositional tautology is true. 

Enayat and Pakhomov in \cite{EnayatPakhomov} proved that actually a very modest fragment of propositional closure, $\PC$, is already enough to yield a non-conservative theory.

\begin{definition} \label{def_dc}
	By \df{Disjunctive Correcntess} ($\DC$), we mean the following principle:
	\begin{displaymath}
	\forall (\phi_i)_{i \leq c} \in \SentSeq_{\LPA}\Big(T \bigvee_{i \leq c} \phi_i \equiv \exists i \leq c T \phi_i \Big).
	\end{displaymath}
\end{definition}
In other words, $\DC$ expresses that any finite disjunction is true iff one of its disjuncts is. Here "finite" is understood in the formalised sense, so that it may refer to nonstandard objects. We treat the symbol $T\bigvee_{i \leq c} \phi_i$ as denoting disjunctions with parentheses grouped to the left for definiteness. 

\begin{theorem}[Enayat--Visser] \label{th_dc_not_conservative}
	$\CT^- + \DC$ is equivalent to $\CT_0$. Consequently, $\CT^-+ \DC$ is not conservative over $\PA$. 
\end{theorem}
This theorem is really striking. Admittedly, $\DC$ can be viewed as a natural extension of compositional axioms. We simply want to allow that the truth predicate behaves compositionally with respect not just to binary (or standard) disjunctions, but to arbitrary finite ones.

\subsection{Disjunctions with stopping conditions}

The main technical tool which we are going to use in this article are disjunctions with stopping conditions, a tool implicitly introduced (but not officially defined), in \cite{smith}. This is a particular propositional construction which is a very useful tool in the analysis of $\CT^-$. The motivation and proofs of the cited facts concerning disjunctions with stopping conditions can be found in \cite{WcisloKossak}.

\begin{definition} \label{def_disjunctions_stopping_cond}
	Let $(\alpha_i)_{i \leq c}, (\beta_i)_{i \leq c}$ be sequences of sentences. We define the \df{disjunction of $\beta_i$ with stopping condition $\alpha$} for $i \in [j,c]$ by backwards induction on $j$:
	\begin{eqnarray*}
	\bigvee_{i = c}^{\alpha,c} \beta_i & = & \alpha_c \wedge \beta_c \\
	\bigvee_{i = j}^{\alpha,c} \beta_i & = & (\alpha_j \wedge \beta_j) \vee (\neg \alpha_j \wedge \bigvee_{i = j+1}^{\alpha,c} \beta_i ).
	\end{eqnarray*} 
\end{definition}
The key feature of disjunctions with stopping conditions is that they allow us to use disjunctive correctness in some very limited range of cases which actually suffice for certain applications without actually committing to the full strength of this axiom.

\begin{theorem} \label{tw_disjunctions_stopping_cond}
	Let $(M,T) \models \CT^-$. Let $(\alpha_i)_{i \leq c}, (\beta_i)_{i \leq c} \in \SentSeq_{\LPA}(M)$ be sequences of sentences. Suppose that $k_0 \in \omega$ is the least number $j$ such that $(M,T) \models T \alpha_{j}$ holds. Then
	\begin{displaymath}
	(M,T) \models T \bigvee_{i = 0}^{\alpha,c} \beta_i \equiv T \beta_{k_0}.
	\end{displaymath}
\end{theorem}
Notice that above we assume that $k_0 \in \omega$, i.e., it is in the standard part of $M$. In other words: if we are guaranteed that some $\alpha_k$ holds for a standard $k$, we can make an infinite case distinction of the form:
"either $\alpha_0$ holds and then $\beta_0$ or $\alpha_1$ holds and then $\beta_1$ ... or $\alpha_c$ holds and then $\beta_c$" so that it actually works correctly in the presence of compositional axioms alone without assuming any induction whatsoever. The proof of Theorem \ref{tw_disjunctions_stopping_cond} (together with applications) may be found in \cite{WcisloKossak}.

The following proposition explains why disjunctions with stopping conditions are so named. 
\begin{proposition} \label{stw_unique_disjunctions_vs_stopping}
	 Suppose that $\alpha_i \beta_i, i \leq c$ are sentences of propositional logic. Then every boolean valuation which makes exactly one of $\alpha_i$ satisfied makes the following equivalence satisfied:
	\begin{displaymath}
	\bigvee_{i=0}^c \alpha_i \wedge \beta_i \equiv \bigvee_{i=0}^{\alpha,c} \beta_i.
	\end{displaymath} 
	Moreover, this is provable in $\PA$. 
\end{proposition}
\begin{proof}
	We work in $\PA$. Fix any valuation which makes exactly one of the sentences $\alpha_i$ true, say, $i=k$.  
	
	It is clear that the disjunction $\bigvee_{i = k}^c \alpha_i \wedge \beta_i$ is equivalent to $\beta_k$. We will show by backwards induction on $j$ that all formulae  $\bigvee_{i = j}^{\alpha_i,c} \beta_i$ are equivalent to $\beta_k$. 
	
	Suppose that $j=k$. Since $\alpha_k$ holds, we immediately have the following equivalence: 
	\begin{displaymath}
	\bigvee_{i=k}^{\alpha,c} \beta_i = (\alpha_k \wedge \beta_k) \vee (\neg \alpha_k \wedge \bigvee_{i=k+1}^{\alpha,c } \beta_i ) \equiv \beta_k.
	\end{displaymath}
	
	Suppose that the claim holds for $j+1 \leq k$. Since $j<k$, by assumption $\alpha_j$ is not true. Hence, again by elementary manipulations, the following equivalence holds:
	 \begin{displaymath}
	 \bigvee_{i=j}^{\alpha,c} \beta_i = (\alpha_j \wedge \beta_j) \vee (\neg \alpha_j \wedge \bigvee_{i=j+1}^{\alpha,c } \beta_i ) \equiv \bigvee_{i=j+1}^{\alpha,c} \beta_i.
	 \end{displaymath}  
	 By induction hypothesis, the last formula is equivalent to $\beta_k$. This proves our claim.
\end{proof}
Theorem \ref{tw_disjunctions_stopping_cond} can be proved by following the above argument, starting with $k_0$ instead of $k$ and noticing that in this case, we only need to perform standardly many steps in of induction, so it can be carried out externally. Let us also remark, that Proposition \ref{stw_unique_disjunctions_vs_stopping} can be clearly proved in much weaker subsystem of $\PA$ such as $\IDelta_0 +\exp$. 

Most importantly for this article, the behaviour of disjunctions with stopping conditions can be partly encoded as a propositional tautology.  

We will use the following notation: if $(\alpha_i)_{i \leq c}$ is a sequence of sentences, then by $\Unique_{i \leq c} \alpha_i$, we mean the following sentence:
\begin{displaymath}
\bigvee_{i \leq c} \Big(\alpha_i \wedge \bigwedge_{j \neq i} \neg \alpha_j \Big).
\end{displaymath}
It expresses that exactly one of $\alpha_i$s is true. 

\begin{cor} \label{cor_disjunctions_stopping_tautologies}
For any sentences $\alpha_i,\beta_i, i \leq c$, the following is a propositional tautology:
	\begin{displaymath}
	\Unique_{i \leq c} \alpha_i \rightarrow \Big(\bigvee_{i = 0}^{c, \alpha} \beta_i \equiv \bigvee_{i=0}^c \alpha_i \wedge \beta_i \Big).
	\end{displaymath}
Morevoer, this is provable in $\PA$. 
\end{cor}

\section{The main result}

In this section, we prove the main result of our paper. We will show that propositional soundness principle added to $\CT^-$ becomes non-conservative (and actually equivalent to $\CT_0$) upon adding an innocuous principle which by itself can be easily shown to be conservative.  

\begin{definition} \label{def_qf_correctness}
	By \df{quantifier-free correctness principle} $(\QFC)$, we mean the following axiom:
	\begin{displaymath}
	\forall \phi \in \qfSent_{\LPA} \Big(T\phi \equiv \Tr_0\phi\Big).
	\end{displaymath}
\end{definition}
In other words, on quantifier-free sentences arithmetical partial truth and truth in the sense of the $T$ predicate agree. Notice that this allows us to use full induction when reasoning about the truth predicate applied to quantifier-free sentences, since the truth predicate restricted to such sentences is equivalent to an arithmetical formula. It turns out that this innocuous principle is enough to yield propositional soundness nonconservative. 

\begin{theorem} \label{tw_propositional_soundness_with_qf_correctness}
	The theory $\CT^- + \QFC + \PS$ is not conservative over $\PA$. In fact, it is exactly equivalent to $\CT_0$. 
\end{theorem}
Crucially, $\CT^- + \QFC$ is by itself conservative over $\PA$. 
\begin{theorem} \label{tw_qfc_conservative}
	The theory $\CT^- + \QFC$ is conservative over $\PA$. 
\end{theorem}
The proof of this fact is a routine application of Enayat--Visser proof of conservativeness of $\CT^-$. For completeness, we present it in Appendix B.

Now, we can present the last crucial ingredient of our proof. As we already mentioned, disjunctive correctness was proved to be equivalent to $\CT_0$ (over $\CT^-$) in \cite{EnayatPakhomov}. However, by inspection of the proof, it can be seen that actually somewhat weaker assumption is employed, as the disjunctive correctness is used only with respect to one rather specific kind of formulae. 

\begin{definition}
	By \df{Atomic Case Distinction Correctness} ($\ACDC$) we mean the following axiom: For any sequence of arithmetical sentences $(\phi_i)_{i \leq c} \in \SentSeq_{\LPA}$ , for any closed term $t \in \ClTerm_{\LPA}$, the following equivalence holds:
	\begin{displaymath}
	T \left(\bigvee_{i \leq c} t= \num{i} \wedge \phi_i\right) \equiv \exists a \leq c \left( \val{t} = a  \wedge T\phi_a \right).
	\end{displaymath} 
\end{definition}

\begin{theorem}[Essentially Enayat--Pakhomov]  \label{tw_atomic_case_distinction_nonconservative}
	$\CT^- + \ACDC$ is equivalent to $\CT_0$. In particular it is not conservative over $\PA$. 
\end{theorem}
As we already mentioned, this theorem is proved by inspection of the earlier argument by Enayat and Pakhomov. For the convenience of the reader, we will discuss it in Appendix A. 

Now, we are ready to present the proof of our main result \ref{tw_propositional_soundness_with_qf_correctness}. 
\begin{proof}[Proof of Theorem \ref{tw_propositional_soundness_with_qf_correctness}]
	Fix any model $(M,T) \models \CT^- + \QFC + \PS.$ We will show that
	\begin{displaymath}
	(M,T) \models \CT^- + \ACDC,
	\end{displaymath}
	which  shows by Theorem \ref{tw_atomic_case_distinction_nonconservative} that $(M,T) \models \CT_0$.
	
	Fix any $c\in M$, a closed term $t \in \ClTerm_{\LPA}(M)$, and an arbitrary sequence of sentences $(\phi_i)_{i\leq c} \in \SentSeq_{\LPA}(M)$. 
	
	First, suppose that there exists $a \leq c$ such that $\val{t} = a$ and $T\phi_a$ holds. Observe that:
	\begin{displaymath}
	(t = a \wedge \phi_a) \rightarrow \bigvee_{i \leq c} t= \num{i} \wedge \phi_i
	\end{displaymath}
	is recognised in $M$ as a propositional tautology. Hence, by $\CT^- + \PS$, we obtain:
	\begin{displaymath}
	(M,T) \models T \left(\bigvee_{i \leq c} t= \num{i} \wedge \phi_i\right).
	\end{displaymath}
	This proves one direction of $\ACDC$. For the harder direction, assume that 
	\begin{displaymath}
	(M,T) \models T \left(\bigvee_{i \leq c} t= \num{i} \wedge \phi_i\right).
	\end{displaymath}
	We first show that indeed $M \models \val{t} \leq c$. Suppose otherwise. Then, this fact is recognised by the partial  arithmetical truth predicate as follows:
	\begin{displaymath}
	M \models \Tr_0 \bigwedge_{i \leq c} \neg t = \num{i}.
	\end{displaymath} 
	By $\QFC$, the same holds for the truth predicate $T$ rather than $\Tr_0$. Moreover, notice that the following sentence is a propositional tautology:
	\begin{displaymath}
	\bigwedge_{i \leq c} \neg t = \num{i} \rightarrow \neg \bigvee_{i \leq c} t= \num{i} \wedge \phi_i.
	\end{displaymath}
	Hence, by propositional soundness $\PS$ and our assumption that $T \bigvee_{i \leq c} t= \num{i} \wedge \phi_i$ holds, the value of $t$, as computed in $M$, is below $c$. 
	
	Now, fix $a \leq c$ such that $\val{t} = a$. Fix any permutation $\sigma: \{0, \ldots, c\} \to \{0, \ldots, c\}$ such that $\sigma(a) = 0$. Since disjunctions are associative and commutative provably in $\PA$ (and in much weaker systems), by propositional soundness $\PS$, the following holds:
	\begin{displaymath}
	(M,T)  \models T \left (\bigvee_{i \leq c} t= \num{i} \wedge \phi_i \right) \equiv T \left (\bigvee_{i \leq c} t= \num{\sigma(i)} \wedge \phi_{\sigma(i)} \right).
	\end{displaymath}
	
	Now, notice that exactly one of the formulae $t = \num{i}$ is true, and this can be expressed as follows:
	\begin{displaymath}
	M \models \Tr_0 \bigvee_{i \leq c} \Bigl( t= \num{i} \wedge \bigwedge_{j \neq i} \neg t = \num{j} \Bigr).
	\end{displaymath}
	By $\QFC$, using our notation from previous section, this is equivalent to:
	\begin{displaymath}
	M \models T \Unique_{i \leq c} t = \num{i}.
	\end{displaymath}
	The same argument applies, if we consider sentences $t = \num{\sigma(i)}$ rather than $ t = \num{i}$. 
	By Corollary \ref{cor_disjunctions_stopping_tautologies}, the following is a propositional tautology, hence true in the sense of the predicate $T$ by $\PS$:
	\begin{displaymath}
	\Unique_{i \leq c} t = \num{i} \rightarrow \left( \left( \bigvee_{i \leq c} t= \num{i} \wedge \phi_i \right) \equiv \bigvee_{i = 0}^{t = \num{i}, c} \phi_i \right).
	\end{displaymath}
	Again, this holds if we consider sequences $t=\num{\sigma(i)}$ and $\phi_{\sigma(i)}$ instead.
	Putting it all together, we know that the following formulae are true:
	\begin{displaymath}
	(M,T) \models T \Unique_{i \leq c} t = \sigma(i) \wedge T \bigvee_{i = 0}^c t= \num{\sigma(i)} \wedge \phi_{\sigma(i)}.
	\end{displaymath}
	Therefore, 
	\begin{displaymath}
	(M,T) \models T \bigvee_{i = 0}^{t= \num{\sigma(i)}, c} \phi_{\sigma(i)}.
	\end{displaymath}
	By Theorem \ref{tw_disjunctions_stopping_cond} on disjunctions with stopping conditions, as the above disjunction stops at $i = 0$, we obtain:
	\begin{displaymath}
	(M,T) \models T \phi_{\sigma(0)}.
	\end{displaymath}
	Since $\sigma(0) = a = \val{t}$, this concludes our argument.
\end{proof}
\section{Appendix A: The strength of $\ACDC$}

In the main part, we crucially used the observation that Atomic Case Distinction Correctness, $\ACDC$ is equivalent to $\CT_0$. As we already mentioned, this result is really due to Enayat and Pakhomov, as this is what their arguments in \cite{EnayatPakhomov} actually show. However, since verifying this claim would be admittedly cumbersome, we will rather repeat their argument below. 

Following closely the presentation in the original paper, we split our argument into two parts. We first show that $\ACDC$ together with internal induction yields $\Delta_0$-induction for the truth predicate. Subsequently, we show that $\ACDC$ implies internal induction. Before any of this happens let us define what internal induction actually is. 

\begin{definition} \label{def_internal_induction}
	By \df{Internal Induction} ($\INT$), we mean the following axiom:
	\begin{displaymath}
	\forall \phi \in \form^{\leq 1}_{\LPA} \Big(\forall x \bigl( T \phi(\num{x}) \rightarrow T \phi(\num{x+1}) \bigr) \rightarrow \forall x \bigl( T\phi(\num{0}) \rightarrow \forall x \ T \phi(\num{x}) \bigr)  \Big).
	\end{displaymath}
\end{definition}
In other words, internal induction expresses that any arithmetical formula satisfies induction under the truth predicate.

\begin{theorem} \label{tw_acdc_plus_int}
	$\CT^- + \ACDC + \INT$ is equivalent to $\CT_0$.
\end{theorem}
\begin{proof}
	It can be directly verified that $\CT_0$ implies $\INT$ and full $\DC$. Therefore, we will focus on the harder direction, showing that $\CT^- + \ACDC + \INT$ implies $\CT_0$.
	
	Fix any model $(M,T) \models \CT^- + \ACDC + \INT$.  We want to show that $(M,T) \models \CT_0$. It is enough to demonstrate that for any $c \in M$, the set $T \cap [0,c]$ is coded, i.e., there exists $s \in M$ such that $a \in T \cap [0,c]$ iff the $a$-th bit of $s$ in the binary expansion is equal to $1$. 
	
	Fix the sequence $(\phi_i)_{i \leq c}$ of sentences such that $\qcr{\phi_i} = i $ if $i$ happens to be an arithmetical sentence (that is, $i \in \Sent_{\LPA}(M)$) and $\phi_i = \qcr{0  \neq 0}$ otherwise. Consider the following formula $\Theta(a,x)$:
	\begin{displaymath}
	\Theta_c(x) := \bigvee_{i \leq c} x = \num{i} \wedge \phi_i. 
	\end{displaymath}
	By $\ACDC$, for $\phi \in \Sent_{\LPA}(M) \cap [0,c]$, 
	\begin{displaymath}
	(M,T) \models T \Theta_c(\num{\phi}) \equiv T \phi.  
	\end{displaymath}
	On the other hand, by $\INT$, the formula $T \Theta_c(\num{x})$ satisfies full induction. In particular, the set of elements smaller than $c$ satisfying this formula is coded. 
\end{proof}

Now, we can move to the second ingredient of the proof:
\begin{theorem} \label{th_acdc_dowodzi_int}
	$\CT^- + \ACDC$ implies $\INT$.
\end{theorem}
In the paper \cite{EnayatPakhomov} which we closely follow in this presentation, the analogue of Theorem \ref{th_acdc_dowodzi_int} is proved by an extremely elegant detour via a theory of iterated truth predicates. 

\begin{definition} \label{def_itb}
	By $\ITB$ (Iterated Truth Biconditionals), we mean a theory with two sorts: a number sort and index sort, over the language with the following symbols:
	\begin{itemize}
		\item The function symbols of $\LPA$, whose arguments come from the number sort. 
		\item A fresh predicate $T(\alpha, x)$, where $\alpha$ comes from the index sort and $x$ from the number sort. We will also denote it with $T_{\alpha}(x)$.
		\item A fresh predicate $\alpha \prec \beta$, whose arguments come from the index sort.
	\end{itemize}
	Its axioms consist of $\PA$, axioms saying that $\prec$ is a linear ordering of the index sort and the following scheme:
	 \begin{displaymath}
	 \forall \alpha \Big(T_{\alpha} \phi \equiv \phi^{\prec \alpha} \Big),
	 \end{displaymath}
	 where $\phi$ comes from the full language and $\phi^{\prec \alpha}$ is $\phi$ with the index-sort quantifiers $\forall \beta, \exists \beta$ replaced with $\forall \beta \prec \alpha, \exists \beta \prec \alpha$. 
\end{definition}

$\ITB$ axiomatises a hierarchy of truth predicates over a linear order. The key point is that this order cannot have infinite descending chains. The theorem below was proved in \cite{EnayatPakhomov}, based on the main result in \cite{visser_yablo}. 
\begin{theorem} \label{tw_itb_plus_nonwf_sprzeczna}
	The theory $\ITB$ together with the axioms $\forall \alpha \exists \beta \ \beta \prec \alpha$ and $\exists \alpha \ \alpha = \alpha$ is inconsistent.
\end{theorem}
By the above theorem, there exists a finite fragment $\Gamma$ of $\ITB$ which proves that $\prec$ has the least element. This theory contains finitely many biconditionals of the form:
 \begin{displaymath}
\forall \alpha \Big(T_{\alpha} \phi \equiv \phi^{\prec \alpha} \Big).
\end{displaymath}
Let $\phi_1,\ldots,\phi_n$ be the enumeration of sentences which occur in the biconditionals from $\Gamma$. Let us denote the biconditional involving $\phi_i$ with $B(\phi_i)$. 

\begin{proof}[Proof of Theorem \ref{th_acdc_dowodzi_int}.]
	Let $(M,T) \models \CT^- + \ACDC.$ Fix any $\phi \in \form_{\LPA}^{\leq 1}(M)$ such that for some $c_0 \in M$, $(M,T) \models  T \phi(\num{c_0})$.  We will show that there exists the least $c \in M$ such that $(M,T) \models T \phi(\num{c})$. Since $\phi$ is arbitrary, and by compositionality of $T$, this implies that internal induction holds in $(M,T)$.
	
	By induction we will construct in $M$ a sequence of interpretations $\iota_a, a \in M$ of $\Gamma \subset \ITB$, i.e., a sequence of tuples of formulae: the definitions of domains for number and index sorts, the interpretations of the arithmetical symbols, and the interpretations for the predicates $\prec$, $T(\alpha,x)$. 
	
	\begin{itemize}
		\item 	For all $a$, $\iota_a$ interprets arithmetical symbols by identity and the domain of number quantifiers is the whole $M$ (i.e., the domain is defined by the formula $x=x$). 
		\item 	The $a$-th domain of index quantifiers is given by $d_a(x) := x \leq a \wedge \phi(x)$.
		\item The index inequality $\prec$ is interpreted by the usual inequality $<$. 
		\item The predicate $T(\alpha,x)$ is defined recursively as follows:
	\end{itemize}
\begin{displaymath}
\bigvee_{i \leq n} \Bigl( x = \phi_i \wedge \bigvee_{j<a} \alpha = \num{j} \wedge \phi(\num{j})  \wedge \iota_j(\phi_i) \Bigr).
\end{displaymath}
	
	We will show that for all $a$, if $(M,T) \models T\phi(\num{a})$, then $\iota_a$ is indeed an interpretation of $\Gamma$  under the truth predicate. This means that for all sentences $\psi \in \Gamma$, 
\begin{displaymath}
(M,T) \models T \iota_a(\psi).
\end{displaymath}

This is immediate for arithmetical axioms and the ordering axioms for $\prec$. Thus it is enough to check that the claim is satisfied for the truth biconditionals. Fix $k \leq n$ and $a \in M$. We want to check that:
\begin{displaymath}
(M,T) \models  T \iota_a \forall \alpha \Big(T_{\alpha} \phi_k \equiv \phi_k^{\prec \alpha} \Big).
\end{displaymath}
	
If there are no $a'<a$ such that $(M,T) \models T\phi(\num{a'})$, then the interpretation of the universal quantifier makes the sentence trivially true. So suppose otherwise and fix any $\alpha< a$ such that $Td_a(\num{\alpha})$ holds. We want to check that the following holds:
\begin{displaymath}
T\iota_a T(\alpha,\phi_k) \equiv T\iota_a \phi_k^{\prec \alpha}.  
\end{displaymath}    
	Expanding the definition of $\iota_a$ on the left-hand side of the equivalence yields:
\begin{displaymath}
T \biggl(\bigvee_{i \leq n} \Bigl( \phi_k = \phi_i \wedge \bigvee_{j<a} \alpha = \num{j} \wedge \phi(\num{j})  \wedge \iota_j(\phi_i) \Bigr) \biggr).
\end{displaymath}
The first disjunction has standardly many disjuncts, of which only one is true, so by compositional axioms, this is equivalent to:
\begin{displaymath}
T \bigvee_{j<a} \alpha = \num{j} \wedge \phi(\num{j})  \wedge \iota_j(\phi_k). 
\end{displaymath}
By $\ACDC$, this is equivalent to:
\begin{displaymath}
T \phi(\num{\alpha}) \wedge T \iota_\alpha(\phi_k).
\end{displaymath}
By assumption on $\alpha$, this is equivalent to 
\begin{displaymath}
T \iota_\alpha(\phi_k).
\end{displaymath}	
Now, it is enough to check that the following equivalence holds:
\begin{displaymath}
T \iota_\alpha(\phi_k) \equiv T \iota_a  \phi_k^{\prec \alpha}.
\end{displaymath}
We essentially check   by induction on complexity of subformulae $\psi$ of $\phi_k$ that this equivalence holds for all $\psi$. To make this more precise, we introduce the following definition. We say that a tuple $t_1, \ldots, t_m \in \ClTerm_{\LPA}(M)$ is \df{suitable} for a formula $\psi$ if $\psi$ has $m$ free variables and for every term $t$ corresponding to an index variable $\beta$, $(M,T) \models T d_{\alpha}(t)$. 

Now, by induction on complexity of formulae, we will show that for any subformula $\psi$ of $\phi_k$, and any suitable tuple $\bar{t}$ of closed terms in the sense of $M$, the following equivalence holds:
\begin{displaymath}
T \iota_{\alpha} (\psi)(\bar{t}) \equiv T \iota_a(\psi)^{\prec \alpha}(\bar{t}).
\end{displaymath}
The induction steps for connectives and number quantifiers, as well as the initial step for the arithmetical atomic formulae and the atomic formula $\beta \prec \gamma$ are immediate. Let us now focus on the initial case for the formula $T(\beta,x)$. Fix any suitable pair of terms $t_1, t_2$. In particular this means that the value of $t_2$ is no greater than $\alpha < a$. $T \iota_\alpha T(t_1,t_2)$ is the following sentence:
\begin{displaymath}
T \biggl(\bigvee_{i \leq n} \Bigl( t_1 = \phi_i \wedge \bigvee_{j< \alpha} t_2 = \num{j} \wedge \phi(\num{j})  \wedge \iota_j(\phi_i) \Bigr) \biggr).
\end{displaymath}
By $\ACDC$ and the fact that $\val{t_2} \leq \alpha < a$, this is equivalent to:
\begin{displaymath}
T \biggl(\bigvee_{i \leq n} \Bigl( t_1 = \phi_i \wedge \bigvee_{j<a} t_2 = \num{j} \wedge \phi(\num{j})  \wedge \iota_j(\phi_i) \Bigr) \biggr).
\end{displaymath} 
(The two formulae differ by the range of the second disjunction.) Since the second formula is equal to $T \iota_a T(t_1,t_2) = T \iota_aT(t_1,t_2)^{\prec \alpha}$, the atomic case is proved. 

What remains to be proved is the induction step for the index quantifier. Suppose that our claim holds for $\psi$ and consider the formula $\forall \beta \psi(\beta).$ Notice that the following equalities hold:
\begin{eqnarray*}
 \iota_\alpha \forall \beta \psi(\beta) & = &  \forall x \Big(d_{\alpha}(x) \rightarrow \iota_{\alpha}\psi(x) \Big) \\
 \iota_a \Big( (\forall \beta \psi(\beta))^{\prec \alpha} \Big) & = &  \forall x \Big(d_{\alpha}(x) \rightarrow \iota_a\psi^{\prec \alpha}(x) \Big).
\end{eqnarray*}
By induction hypothesis and the compositional axioms if we substitute suitable terms in the formulae on the right-hand side, then the first one is true if and only if the second one is. This concludes the induction argument and the whole proof. 
\end{proof}

\section*{Appendix B: Conservativeness of $\CT^- + \QFC$}

In the main part, we claimed that the quantifier-free correctness can be added to $\CT^-$ still yielding a conservative theory. As we already noted, this is a very simple application of the Enayat--Visser construction, but we could not find this exact statement in the literature.\footnote{Similar statements concerning satisfaction classes containing $\Sigma_n$ arithmetical truth can be found e.g. in \cite{engstrom}, but the definitions of satisfaction class there is slightly different from the one we use. However, the conservativeness result proved here is neither surprising nor really original.} Therefore, we decided to include a proof of our statement. However, the reader should feel entirely free to skip it. 

\begin{definition} \label{def_partial_truth_predicate}
		Let $M \models \PA$. We say that a set $T_0 \subset \Sent_{\LPA}(M)$ is a \df{partial compositional truth predicate} if the following conditions hold:
	\begin{itemize}
		\item For any $\phi \in T_0$ and any $\psi$ which results by substituting closed terms into a direct subformula of $\phi$, $\psi \in T_0$. 
		\item If $\psi \in T_0$, then  the sentences which result by substituting closed terms into direct subformulae of $\phi$ satisfy compositional axioms 1-4. of $\CT^-$.  
		\item $T_0$ satisfies extensionality axiom 5. of $\CT^-$. 
	\end{itemize}
\end{definition}

We will derive Theorem \ref{tw_qfc_conservative} from the following, more general fact. 

\begin{theorem} \label{th_extensions_of_truth_predicates}
 Let $M_0 \models \PA$ and let $T_0 \subset M_0$ be a partial truth predicate. Then there exists an elementary extension $(M_0,T_0) \preceq (M',T)$ and $T' \supseteq T$ such that $(M',T') \models \CT^-$. 
\end{theorem} 

\begin{proof}[Proof of Theorem \ref{tw_qfc_conservative} from Theorem \ref{th_extensions_of_truth_predicates}]
	Let $M \models \PA$ and let $T_0 \subset M_0$ be defined as the set os sentences $\phi$ such that $M \models \Tr_0(\phi).$ We apply Theorem \ref{th_extensions_of_truth_predicates} to $(M_0,T_0)$ obtaining am elementary extension $(M',T) \succeq (M_0,T_0)$ and  $T' \supseteq T$ such that $(M',T') \models \CT^-$.  
	
	Now, observe that actually $(M',T') \models \CT^- +  \QFC$. Indeed, by elementarity $T$ is exactly the set of $\phi \in \Sent_{\LPA}(M')$ such that $M' \models \Tr_0(\phi)$. 
\end{proof}

Now we turn to the proof of Theorem \ref{th_extensions_of_truth_predicates}. Since we are dealing with truth predicates for a language with terms and we include extensionality in our axioms, we have to take care of certain additional technicalities. Before we proceed to the proof, we will introduce some definitions and notation. 

\begin{definition} \label{def_template}
	Let $M \models \PA$ and let $\phi \in \form_{\LPA}(M)$. By a \df{trivialisation} of $\phi$, we mean a formula $\widehat{\phi}$ such that:
	\begin{itemize}
		\item There exists a sequence of terms $\bar{t} \in \TermSeq_{\LPA}(M)$ such that $\phi = \widehat{\phi}(\bar{t})$. 
		\item No variable occurs in $\widehat{\phi}$ both free and bound.
		\item No free variable occurs in $\widehat{\phi}$ more than once.
		\item No closed term occurs in $\widehat{\phi}$. 
		\item No complex term whose all variables are free occurs in $\widehat{\phi}$.
		\item $\widehat{\phi}$ is the least formula with the above properties. (In order to guarantee uniqueness.) 
	\end{itemize} 
\end{definition}

For instance, if $\phi = \exists x \forall y \Big( x +( z \times S0 + 0 \times u ) = x \times y + 0 \Big)$, then 
\begin{displaymath}
\widehat{\phi} = \exists x \forall y \Big( x + v_1 = x \times y + v_2 \Big),
\end{displaymath}
where $v_1, v_2$ are chosen so as to minimise the formula $\widehat{\phi}$. 

\begin{itemize}
	\item We say that two formulae $\phi_1, \phi_2$ are \df{syntactically similar} if $\widehat{\phi_1} = \widehat{\phi_2}.$ We denote it with $\phi_1 \sim \phi_2$.
	\item If $\phi \in \form_{\LPA}$ and $\alpha \in \Asn(\phi)$, then by $\phi[\alpha]$ we mean the sentence resulting by substituting the numeral $\num{\alpha(v)}$ for each variable $v$. 
	\item If $\phi_1, \phi_2 \in \form_{\LPA}, \alpha_1 \in \Asn(\phi_1), \alpha_2 \in \Asn(\phi_2)$, then we say that $(\phi_1,\alpha_1)$ is \df{extensionally equivalent} to $(\phi_2, \alpha_2)$ if $\phi_1 \sim \phi_2$ and there exist two sequences of closed terms $\bar{t_1}, \bar{t_2} \in \ClTermSeq_{\LPA}$ such that $\overline{\val{t_1}} = \overline{\val{t_2}}$ (the values of terms in $\bar{t_1}, \bar{t_2}$ are pointwise equal), $\phi_1 = \phi(\bar{t_1}), \phi_2 = \phi(\bar{t_2})$, where $\phi = \widehat{\phi_1}= \widehat{\phi_2}$. We denote this relation with $(\phi_1, \alpha_1) \sim (\phi_2, \alpha_2)$. 
\end{itemize}
 
Notice that the syntactic similarity and extensional equivalence are both equivalence relations.

\begin{proof}
	Let $M_0$ be any model of $\PA$ and let $T_0 \subset M_0$ be a partial truth predicate. We will construct a chain of models $(M_i,T_i,S_i), i \in \omega$. The chain of models $(M_i,T_i)$ will be elementary and the binary predicate $S_i$ will be partial satisfaction predicates extending one another and extending $T_i$.

	We perform the construction in the following way: once we have constructed the model $(M_i,T_i,S_i)$, we let $(M_{i+1},T_{i+1}, S_{i+1})$ be any model of the theory $\Theta_{i+1}$ consisting of the following axioms in the arithmetical language with additional predicates $S_{i+1}, T_{i+1}$:
	\begin{itemize}
		\item $\ElDiag(M_i,T_i)$. (The elementary diagram of $(M_i, T_i)$, formulated with $T_{i+1}$ replacing $T_i$.)
		\item $\Comp(\phi), \phi \in \form_{\LPA}(M_i)$. (The compositionality scheme, to be defined later.) 
		\item $\forall \phi,\phi' \in \form_{\LPA} \forall \alpha \in \Asn(\phi), \alpha' \in \Asn(\phi') \ \Big(S_{i+1}(\phi,\alpha) \equiv S_{i+1}(\phi',\alpha') \Big)$. (The extensionality axiom)
		\item $\forall x \Big(T_{i+1}(x) \rightarrow S_{i+1}(x,\emptyset)\Big)$. (The satisfaction predicate $S_{i+1}$ agrees with $T_{i+1}$.)
		\item $S_{i+1}(\phi,\alpha)$, where $\phi \in \form_{\LPA}(M_{i-1})$, $\alpha \in \Asn(\phi)$ and $(\phi,\alpha) \in S_i$. (The preservation scheme.)
	\end{itemize}

An instance of the compositionality scheme $\Comp(\phi)$ is defined as the disjunction of the following clauses:

\begin{enumerate}
	\item $\exists s, t \in \Term_{\LPA} \Big( \phi = (s=t) \wedge \forall \alpha \in \Asn(\phi) \  S(\phi, \alpha) \equiv s^{\alpha} = t^{\alpha}\Big).$
	\item $\exists \psi \in \form_{\LPA} \Big( \phi = (\neg \psi) \wedge \forall \alpha \in \Asn(\phi) \ S(\phi, \alpha) \equiv \neg S(\psi,\alpha) \Big).$
	\item $\exists \psi, \eta \in \form_{\LPA} \Big(\phi = (\psi \vee \eta) \wedge \forall \alpha \in \Asn(\phi) \ S(\phi,\alpha) \equiv S(\psi, \alpha) \vee S(\eta, \alpha)  \Big).$
	\item $\exists \psi \in \form_{\LPA}, v \in \Var \Big(\phi = (\exists v \psi) \wedge \forall \alpha \in \Asn(\phi) \Big(S(\phi,\alpha) \equiv \exists \beta \sim_{v} \alpha \ S(\psi,\beta)\Big)\Big).$
\end{enumerate}

For the time being, suppose that all theories $\Theta_n$ are consistent. We will finish the proof under this assumption and return to it afterwards.

Let $M' = \bigcup_{n \in \omega} M_n$, $T = \bigcup T_n$. Let 
\begin{displaymath}
T' = \set{\phi \in \Sent_{\LPA}(M)}{\exists n \in \omega \ \phi \in M_n \wedge (\phi, \emptyset) \in S_{n+1}}.
\end{displaymath} 
It can be directly verified that the predicate $T'$ defined in such a way satisfies axioms of $\CT^-$ thanks to the assumption that the predicates $S_{n}$ satisfy the compositionality scheme together with preservation and extensionality axioms. Similarly, we check that $T' \supset T$, because each of the predicates $S_n$ extends $T_n$. The details are rather straightforward. The reader can consult the Appendices in \cite{WcisloKossak} or \cite{loccoll}, where a very similar construction is presented.

We have yet to check by induction that all theories $\Theta_n$ are consistent. So assume that this is true for a given $\Theta_n$ and let $M_n \models \Theta_n$. In order to make the proof work uniformly for the successor and the initial steps of induction, we set by convention $M_{-1} = T_{-1} = S_{-1} = \emptyset$. 

We will prove consistency of $\Theta_{n+1}$ in the following way. Consider any finite subtheory $\Gamma \subset \Theta_n$. In the model $M_n$, we will find a binary relation $S$ which satisfies $\Gamma$. 
	
Since $\Gamma$ is finite, there are only finitely many formulae $\phi_1, \ldots, \phi_k$ which occur in the compositionality scheme. 

Consider the equivalence classes $[\phi_i]$ of the formulae $\phi_i$ under the similarity relation $\sim$. Let $\unlhd$ be the transitive closure of the following relation on classes: $[\phi] \unlhd [\psi]$ if there exist $\phi' \in [\phi], \psi' \in [\psi]$ such that $\phi'$ is a direct subformula of $\psi'$. This is indeed an ordering: transitivity and reflexivity is clear, so it is enough to check weak antisymmetry: however this is clear, since if $[\phi] \unlhd [\psi]$, then the total number of connectives and quantifiers in $\phi$ is no greater than in $\psi$. 

We define the extension of $S$ as follows. We define the set $S^0$ by the following conditions. A pair $(\phi,\alpha)$ belongs to $S^0$ if one of the following conditions is satisfied:
\begin{itemize}
	\item $[\phi] \cap M_{n -1} \neq \emptyset$ and $(\phi',\alpha') \sim (\phi,\alpha)$ for some $\phi' \in M_{n-1}$ and $\alpha' \in M_n$ such that $(\phi',\alpha') \in S_n$. 
	\item There exists $\phi' \in M_n$ such that $(\phi',\emptyset) \sim (\phi,\alpha)$ and $\phi' \in T_n$.  
	\item $\phi$ is an atomic formula of the form $t=s$ for some terms $t=s$ and $t^{\alpha} = s^{\alpha}.$ 
\end{itemize}
In the above list, we do not explicitly include the case when $[\phi]$ is minimal among $[\phi_i]$ with respect to the relation $\unlhd$ and $[\phi] \cap M_{n-1} = \emptyset$, but we also implicitly treat this case as covered. Such formulae are simply not satisfied under any assignment. Hence, they effectively define the empty set under the satisfaction predicate. 

Then we inductively construct a series of predicates $S^j$. We define $S^{j+1}$ as the union of $S^j$ with the set of $(\phi,\alpha)$ such that $[\phi] = [\phi_i]$ for some $i \leq k$, $[\phi_i]$ is not minimal with respect to the relation $\unlhd$, and $\phi$ satisfies one of the following conditions: 
\begin{itemize}
	\item There exists $\psi \in M_n$ such that $\phi = \neg \psi$ and $(\psi,\alpha) \notin S^j$.
	\item There exist $\psi, \eta \in M_n$ such that $\phi = \psi \vee \eta$ and $(\psi,\alpha) \in S^j$ or $(\eta,\alpha) \in S^j$.  
	\item There exists $\psi,v \in M_n$ such that $\phi = \exists v \psi$, and $\beta \sim_v \alpha$ such that $(\psi,\beta) \in S^j$.
\end{itemize}

Since we considered only finitely many classes $[\phi_i]$, the construction terminates at some point. Let $S$ be the predicate obtained as the final one in this construction. We claim that $(M_n,T_n,S)$ satisfies $\Gamma$. The elementary diagram of $(M_n,T_n)$ is obviously satisfied in the obtained model. Our construction and the fact that $S_n$ and $T_n$ were compositional and extensional immediately guarantee that the constructed predicate $S$ agrees with $T_n$, preserves $S_n$ for formulae from $M_{n-1}$ and satisfies the instances of the compositional  scheme from $\Gamma$. Finally, we check by induction on $j$ that each $S^j$ satisfies the extensionality axiom. This concludes the proof of consistency of $\Gamma$, the proof of consistency of $\Theta_{n+1}$ and consequently, the proof of Theorem \ref{th_extensions_of_truth_predicates}.	
\end{proof}

\section*{Acknowledgements}
We are grateful to Ali Enayat for a number of helpful comments.  This research was supported by an NCN MAESTRO grant 2019/34/A/HS1/00399 "Epistemic and Semantic Commitments of Foundational Theories."

\end{document}